\documentstyle[12pt]{article}

\newcommand{\be}{\begin{equation}}     
\newcommand{\ee}{\end{equation}}

\newcommand{\bea}{\begin{eqnarray}}
\newcommand{\eea}{\end{eqnarray}}

\newcommand{\eqq}{\!\! & = & \!\!}

\newcommand{\hst}[1]{\rule{#1}{0mm}}

\newcommand{\reff}[1]{(\ref{#1})}    

\newcommand{\vl}[2]{\rule[#1]{0.1mm}{#2}}

\newcommand{\Co}{\mbox{C\hst{-0.44em}\vl{0.1ex}{1.4ex}\hst{0.44em}}}

\newcommand{\comp}{\!\stackrel{\textstyle{\hst{0ex} \atop \circ}}{\hst{0ex}}\!}

\newcommand{\Ga}{\Gamma}
\newcommand{\De}{\Delta}
\newcommand{\al}{\alpha}
\newcommand{\bet}{\beta}
\newcommand{\ga}{\gamma}
\newcommand{\de}{\delta}
\newcommand{\ve}{\varepsilon}
\newcommand{\vr}{\varrho}
\newcommand{\si}{\sigma}
\newcommand{\vp}{\varphi}
\newcommand{\om}{\omega}
\newcommand{\Gas}{\Gamma^{\ast}}
\newcommand{\CGS}{\Co\Gas}
\newcommand{\CG}{\Co\Ga}
\newcommand{\la}{\langle}
\newcommand{\ra}{\rangle}
\newcommand{\rar}{\rightarrow}
\newcommand{\lra}{\longrightarrow}
\newcommand{\ot}{\otimes}

\begin{document}

\begin{titlepage}

\hspace*{\fill} \parbox{8em}{BONN--TH--95--16 \\
                             August 1995        }

\vspace{18mm}

\begin{center}
{\LARGE\bf Universal $R$--matrices \\[-1.5ex]
           for finite Abelian groups --- \\[1ex]
           a new look at graded multilinear algebra} \\
\vspace{11mm}
{\large M. Scheunert} \\
\vspace{5mm}
Physikalisches Institut der Universit\"{a}t Bonn \\
Nu{\ss}allee 12, D--53115 Bonn, Germany
\end{center}

\vspace{20mm}
\begin{abstract}
\noindent
The universal $R$--matrices and, dually, the coquasitriangular structures
of the group Hopf algebra of a finite Abelian group (resp. of an arbitrary
Abelian group) are determined. This is used to formulate graded multilinear
algebra in terms of triangular or cotriangular Hopf algebras. For the
convenience of the reader, in a separate section the definitions and basic
properties of quasitriangular and coquasitriangular Hopf algebras are
recalled.
\end{abstract}

\vspace{\fill}

q-alg/9508016

\end{titlepage}

\setcounter{page}{2}

\section{Introduction}
In the present work I am going to discuss the universal $R$--matrices [1]
for finite Abelian groups, or, more precisely, the universal $R$--matrices
for the group Hopf algebras of these groups [2]. To explain how I became
interested in this problem I remind the reader of the graded multilinear
algebra described in Ref. [3]. As shown in this reference, given an
(arbitrary) Abelian group $\Ga$ and a commutation factor $\si$ on $\Ga$,
the basic constructions of classical multilinear algebra can easily be
generalized to the $\Ga$--graded setting. In modern terms, this means that
the $\Ga$--graded vector spaces form a tensor category (depending on $\si$)
with an internal hom, and this category is rigid if we consider only
finite--dimensional vector spaces. It is this category where the so--called
Lie colour algebras (also called $\si$ Lie algebras) live [4,\,5], and that
is why I have investigated these structures.

Since the present discussion is only meant to serve as a motivation, I do not
want to explain the various category theoretical terms used above, but refer
to Refs. [6,\,7,\,8], where more details are given. However, I would like to
emphasize that the properties of the categories in question are closely
related to the axioms for a triangular Hopf algebra [1]. In fact, given any
triangular Hopf algebra $H$, the (left) $H$--modules form a category of this
type, and conversely, under suitable assumptions, such a category is
equivalent to the category of $H$--modules for some triangular Hopf algebra
$H$ (this is the content of various so--called reconstruction theorems).

In the simple case under consideration we do not want to invoke one of the
reconstruction theorems in order to find the appropriate triangular Hopf
algebra. Instead, we recall that, for a finite Abelian group $\Ga$, a
$\Ga$--graded vector space is nothing but a vector space endowed with a
representation of $\Gas$, the character group of $\Ga$ (see Section 2). Thus
the group Hopf algebra $\CGS$ of $\Gas$ is the natural candidate for $H$,
and we have to determine the triangular universal $R$--matrices for $\CGS$.
(Recall that, for a finite Abelian group $\Ga$, the groups $\Ga$ and $\Gas$
are isomorphic.) Actually, we do more and determine all the universal
$R$--matrices for $\CGS$. From the category theoretical point of view this
means the transition from (symmetric) tensor categories to quasitensor or
braided monoidal categories. The outcome is completely satisfactory: There
exists a simple bijection of the set of bicharacters on $\Ga$ onto the set of
universal $R$--matrices for $\CGS$, and this bijection maps the commutation
factors onto the triangular universal $R$--matrices.

The proof given in Section 2 is elementary, it only uses some simple
properties of finite Abelian groups and their characters. Nevertheless,
it is worth--while to consider the problem from a dual point of view. It is
easy to see that the group Hopf algebras $\CG$ and $\CGS$ are dual to each
other. Thus the $\Ga$--graded vector spaces, i.e., the $\CGS$--modules, can
be identified with the $\CG$--comodules, and the universal $R$--matrices for
$\CGS$ correspond to coquasitriangular structures on $\CG$. Once this has
been realized, the results summarized above follow immediately. Moreover,
in the dual language, these results hold also for not necessarily finite
Abelian groups. All this will be shown in Section 4, after some of the
notions used there (and above) have been explained in the preparatory
Section 3. The final Section 5 contains some concluding remarks and draws
the reader's attention to related work by other authors.

In closing I would like to mention some notations and conventions which will
be used throughout this work. The base field will always be the field $\Co$
of complex numbers, and the multiplicative group of non--zero complex
numbers will be denoted by $\Co_{\ast}$\,. Moreover, $\Ga$ will be an Abelian
group which, for convenience, will be written {\it multiplicatively}. A
character of $\Ga$ is a homomorphism of $\Ga$ into $\Co_{\ast}$\,. The
characters of $\Ga$ form a (multiplicative) group $\Gas$, and the canonical
pairing of $\Gas$ and $\Ga$ will be denoted by a pointed bracket:
\be  \la\,\;,\;\ra : \Gas \times \Ga \lra \Co_{\ast} \ee
\be  \la\ga',\ga\ra = \ga'(\ga)
                      \;\;\mbox{for all}\;\; \ga'\in\Gas\,,\;\ga\in\Ga. \ee
Of course, if $\Ga$ is finite, its characters take their values in the unit
circle $U(1)$.

\section{Universal $R$--matrices for the group
    \newline
         Hopf algebra of a finite Abelian group}
\setcounter{equation}{0}

{\it In the present section the group} $\Ga$ {\it is assumed to be finite.}

Let us first show that there exists a bijective correspondence between
$\Ga$--graded vector spaces and $\Gas$--modules, i.e., vector spaces endowed
with a representation of $\Gas$. In fact, suppose $V$ is a $\Gas$--module
and let $\ga'\cdot\,x$ denote the action of an element $\ga'\in\Gas$ on an
element $x\in V$. Since $\Gas$ is finite, any representation of $\Gas$ is
completely reducible, and since $\Gas$ is Abelian, its irreducible
representations are one--dimensional, hence given by its characters, i.e.,
by the elements of $\Ga$. This implies that $V$ has the decomposition
\be      V = \bigoplus_{\ga\in\Ga}V_{\ga}\:,  \label{dec}   \ee
where the subspaces $V_{\ga}$ of $V$ are defined by
\be  V_{\ga} = \{x \in V\,|\,\ga'\cdot\,x = \la\ga',\ga\ra\,x
                                \;\;\mbox{for all}\;\;\ga'\in\Gas\}\,. \ee
Conversely, suppose that $V$ is a $\Ga$--graded vector space, i.e., suppose
that a decomposition like \reff{dec} is given. Then the action of $\Gas$ on
$V_{\ga}$ defined by
\be  \ga'\cdot\,x = \la\ga',\ga\ra\,x\;\;\mbox{for all}\;\;\ga'\in\Gas\,,\;
                                                            x\in V_{\ga} \ee
makes each of the subspaces $V_{\ga}$ of $V$ into a $\Gas$--module, and
$V$ can be considered as the direct sum of these modules. It is easy to see
that these transitions from $\Gas$--modules to $\Ga$--graded vector spaces
and vice versa are inverse to each other. Moreover, the homomorphisms of
$\Gas$--modules are exactly the homomorphisms of $\Ga$--graded vector spaces,
i.e., the linear mappings which are homogeneous of degree zero.

Now let $\CGS$ be the group algebra of $\Gas$. Recall that $\CGS$ is the
associative algebra consisting of the formal linear combinations of elements
of $\Gas$, with the multiplication inherited from $\Gas$, recall also that
$\Gas$--modules and $\CGS$--modules are essentially the same thing. It is
well--known [2] that $\CGS$ is a Hopf algebra, where the coproduct $\Delta$,
the counit $\varepsilon$, and the antipode $S$ are given (on the basis
$\Gas$ of $\CGS$) by
\bea       \De(\ga') \eqq \ga'\ot\ga' \\
           \ve(\ga') \eqq 1 \\
             S(\ga') \eqq \ga'^{-1}     \eea
for all $\ga'\in\Gas$.

Let us now proceed to the determination of the universal $R$--matrices for
$\CGS$. To find a good Ansatz, let $\si$ be a commutation factor on $\Ga$
and let $V$ and $W$ be two $\Ga$--graded vector spaces. Then, according to
Ref. [3], the canonical isomorphism
\be  \psi_{V,W} : V \ot W \lra W \ot V    \ee
is given by
\be  \psi_{V,W}(x \ot y) = \si(\al,\bet)\,y \ot x \label{col}  \ee
for all $x \in V_{\al}$\,, $y \in W_{\bet}$\,; $\al,\bet \in \Ga$.

On the other hand, for a triangular Hopf algebra $H$ with universal
$R$--matrix
\be  R = \sum_{i}R^1_i \ot R^2_i  \ee
(with $R^1_i, R^2_i \in H$) and for any two $H$--modules $V$ and $W$, the
corresponding braiding isomorphism is fixed by
\be  \psi_{V,W}(x \ot y) = \sum_{i}R^2_i\cdot y \ot R^1_i\cdot x
                                                           \label{braid} \ee
for all $x \in V$, $y \in W$, where, as before, the dot denotes the action of
an element of $H$ on some element of an $H$--module.

In our case, we have $H = \CGS$, and a general element of $\CGS \ot \CGS$
can be uniquely written in the form
\be  R = \sum_{\al',\,\bet'\in \Gas}c_{\al',\,\bet'}\,\al'\ot \bet' \; ,  \ee
with some coefficients $c_{\al',\,\bet'} \in \Co$\,. Inserting this expression
into Eq. \reff{braid} and equating the right hand sides of Eqs. \reff{col}
and \reff{braid} we obtain the condition
\be  \si(\al,\bet) = \sum_{\al',\,\bet'}c_{\al',\,\bet'}\, \la\al',\al\ra
                                             \la\bet',\bet\ra \label{cond} \ee
for all $\al,\bet \in \Ga$.

Using the well--known orthogonality and completeness relations for the
characters, i.e.,
\be  \sum_{\al \in \Ga}\la\al',\al\ra\la\bet',\al\ra^{\ast}
                                  = n \de_{\al',\,\bet'} \label{orth}  \ee
\be  \sum_{\al' \in \Gas}\la\al',\al\ra\la\al',\bet\ra^{\ast}
                                  = n \de_{\al,\,\bet}\;, \label{compl} \ee
where $\al',\bet' \in \Gas$ and $\al,\bet \in \Ga$ and where $n$ denotes the
order of $\Ga$, we can solve Eq. \reff{cond} for the coefficients
$c_{\al',\,\bet'}$ and obtain the following Ansatz for the universal
$R$--matrix
correseponding to $\si$\,:
\be R_{\si} = \frac{1}{n^2}\sum_{\al,\bet \in \Ga \atop \al',\,\bet' \in \Gas}
               \si(\al,\bet)\la\al',\al\ra^{\ast}\la\bet',\bet\ra^{\ast}\,
                                            \al' \ot \bet'\,. \label{ans} \ee

It turns out that Eq. \reff{ans} gives a nice parametrization of the elements
of $\CGS \ot \CGS$. Let ${\cal F}(\Ga \times \Ga)$ be the algebra of all
complex--valued functions on $\Ga \times \Ga$. Then we can prove the
following proposition. \\

\noindent
{\bf Proposition 2.1} \\
{\it 1) The mapping $\si \rar R_{\si}$ of ${\cal F}(\Ga \times \Ga)$ into}
$\CGS \ot \CGS$ {\it given by Eq. \reff{ans} is an algebra isomorphism.}

\noindent
In the subsequent statements, $\si$ denotes an arbitrary element of
${\cal F}(\Ga \times \Ga)$. \\
{\it 2) We have
\be  (\De \ot id)(R_{\si}) = (R_{\si})_{13}(R_{\si})_{23} \label{UR1} \ee
if and only if
\be  \si(\al\bet,\ga) = \si(\al,\ga)\si(\bet,\ga)\;\;\mbox{\it for all}\;\;
                                   \al,\bet,\ga \in \Ga\,. \label{bich1} \ee
3) We have
\be  (id \ot \De)(R_{\si}) = (R_{\si})_{13}(R_{\si})_{12} \label{UR2} \ee
if and only if
\be  \si(\al,\bet\ga) = \si(\al,\bet)\si(\al,\ga)\;\;\mbox{\it for all}\;\;
                                   \al,\bet,\ga \in \Ga\,. \label{bich2} \ee
4) Let} $T:\CGS \ot \CGS \lra \CGS \ot \CGS$ {\it be the usual twist mapping
and define the function $\si^{T} \in {\cal F}(\Ga \times \Ga)$ by
\be  \si^{T}(\al,\bet) = \si(\bet,\al)\;\;\mbox{\it for all}\;\;
                                                   \al,\bet \in \Ga\,. \ee
Then we have
\be   T(R_{\si}) = R_{\si^{T}}\,.  \ee}

Proof \\
By direct calculations, using Eqs. \reff{orth} and \reff{compl}. (The
notation $(R_{\si})_{12}$ etc. is standard, see Section 3.) \\

Finally, we note that
\be  R_{\si}\Delta(h) = (T\Delta(h))R_{\si} \label{UR3} \ee
for all $h \in \CGS$ and all $\si \in {\cal F}(\Ga \times \Ga)$. This is
obvious since $\CGS$ is commutative and cocommutative.

The foregoing results immediately imply the assertion made in the
introduction. \\

\noindent
{\bf Corollary 2.1} \\
{\it Let $\si$ be an arbitrary complex--valued function on $\Ga \times \Ga$.
Then $R_{\si}$ is a universal $R$--matrix for} $\CGS$ {\it if and only if
$\si$ is a bicharacter on $\Ga$. Any universal $R$--matrix for} $\CGS$
{\it can uniquely be written in this form, and $R_{\si}$ is triangular if
and only if $\si$ is skew--symmetric in the sense that
\be  \si(\al,\bet)\si(\bet,\al) = 1\;\;\mbox{for all}
                                   \;\;\al,\bet \in \Ga\,. \label{skew} \ee
The braiding isomorphisms corresponding to a bicharacter $\si$ and to the
associated universal $R$--matrix $R_{\si}$ coincide.}

Proof \\
By definition, $R_{\si}$ is a universal $R$--matrix if it has an inverse and
satisfies the Eqs. \reff{UR1}, \reff{UR2}, and \reff{UR3}. This is the case
if and only if $\si$ takes its values in $\Co_{\ast}$ and satisfies the Eqs.
\reff{bich1} and \reff{bich2}, i.e., if and only if $\si$ is a bicharacter
on $\Ga$. Moreover, $R_{\si}$ is triangular if and only if
\be   T(R_{\si}) = R^{-1}_{\si}\,, \ee
i.e., if and only if
\be   \si^{T} = \si^{-1}\,. \ee
This is exactly Eq. \reff{skew}. The last statement of the corollary has
been used to find the correct Ansatz for $R_{\si}$\,. \\

{}From a calculational point of view, the parametrization \reff{ans} of the
universal $R$--matrices of $\CGS$ is not quite convenient, for it involves
the detour over $\Ga$ and hence a fourfold sum. This drawback can easily be
avoided, as follows. \\

\noindent
{\bf Proposition 2.2} \\
{\it Let $\De_1$ and $\De_2$ be two subgroups of $\Gas$ such that there
exists a non--degenerate pairing}
\be   \tau : \De_1 \times \De_2 \lra \Co_{\ast}\,. \ee
{\it (Since $\Gas$ is finite this is the case if and only if $\De_1$ and
$\De_2$
are isomorphic.) Let $m$ be the common order of $\De_1$ and $\De_2$\,. Then
the equation
\be  \si(\al,\bet) = \frac{1}{m}\sum_{\al_i \in \De_i}\tau(\al_1,\al_2)
                                        \la\al_1,\al\ra\la\al_2,\bet\ra  \ee
defines a bicharacter $\si$ on $\Ga$ and
\be  R^{\tau} = \frac{1}{m}\sum_{\al_i \in \De_i}\tau(\al_1,\al_2)\,
                                            \al_1 \ot \al_2 \label{Rtau} \ee
is the corresponding universal $R$--matrix $R_{\si}$ for} $\CGS$. {\it Every
bicharacter on $\Ga$ and hence every universal $R$--matrix for} $\CGS$ {\it
can uniquely be obtained in this way.}

Proof \\
Needless to say, a non--degenerate pairing $\tau$ of two Abelian groups
$\De_1$ and $\De_2$ is a bicharacter
\be   \tau : \De_1 \times \De_2 \lra \Co_{\ast}  \ee
such that, for any $\al_1 \in \De_1$\,, the condition
\be  \tau(\al_1,\al_2) = 1\;\;\mbox{for all}\;\;\al_2 \in \De_2  \ee
implies that $\al_1$ is the unit element of $\De_1$\,, and likewise with the
roles of $\De_1$ and $\De_2$ interchanged.

The proof of the proposition is quite straightforward, it follows from the
standard duality theory for finite Abelian groups. Thus we shall only give
a few hints.

Let $\si$ be a bicharacter on $\Ga$. Obviously, our task is to investigate
the function on $\Gas \times \Gas$ defined by
\be  (\al',\bet') \lra \sum_{\al,\bet \in \Ga}\si(\al,\bet)
                        \la\al',\al\ra^{\ast}\la\bet',\bet\ra^{\ast}\,. \ee
Define the (left and right) kernels of $\si$ by
\bea  N_1 \eqq \{\al \in \Ga\,|\, \si(\al,\bet) = 1\;\;\mbox{for all}\;\;
                                                       \bet \in \Ga \}  \\
      N_2 \eqq \{\bet \in \Ga\,|\, \si(\al,\bet) = 1\;\;\mbox{for all}\;\;
                                                      \al \in \Ga \}\,. \eea
Then $\si$ induces a non--degenerate pairing $\tilde{\si}$ of $\Ga/N_1$ and
$\Ga/N_2$\,, hence the orders of these groups coincide. We denote them by $m$.
Define the subgroups $\De_1$ and $\De_2$ of $\Gas$ by
\be  \De_i = \{\ga' \in \Gas\,|\, \la\ga',\ga_i\ra = 1\;\;\mbox{for all}\;\;
                                                  \ga_i \in N_i \}\,. \ee
It is well--known that $\De_i$ can be identified with the character group of
$\Ga/N_i$\,, thus $\tilde{\si}$ induces a non--degenerate pairing $\tau$ of
$\De_1$ and $\De_2$\,. This is the one we are looking for.

To be more explicit, define the complex--valued function $\si'$ on
$\Gas \times \Gas$ by
\be  \si'(\al',\bet') = \frac{m}{n^2}\sum_{\al,\bet \in \Ga}\si(\al,\bet)
                        \la\al',\al\ra^{\ast}\la\bet',\bet\ra^{\ast}\,, \ee
where the normalization factor
\be  \frac{m}{n^2} = \frac{1}{m}\frac{1}{(n/m)^2}  \ee
has been chosen conveniently. Using the duality theory for finite Abelian
groups it can be shown that $\si'(\al',\bet')$ vanishes whenever
$\al' \notin \De_1$ or $\bet' \notin \De_2$\,. Let $\tau$ be the restriction
of $\si'$ onto $\De_1 \times \De_2$\,. Then $\tau$ is the pairing described
abstractly above and has all the properties stated in the proposition. \\

\noindent
{\bf Remark} \\
The universal $R$--matrices given by Eq. \reff{Rtau} with
$\De_1 = \De_2 = \Gas$ and $\tau$ a non--degenerate commutation factor on
$\Gas$ have also been obtained in Ref. [9]. \\

\noindent
{\bf Example} \\
Let $n \ge 1$ be an integer, let $\Ga_n$ be the (multiplicative) group of
all $n$th roots of unity, and let $Z_n$ be the ring of integers modulo $n$\,.
Obviously, for any $\al \in \Ga_n$ and $r \in Z_n$\,, the power $\al^r$ has
a well--defined meaning. Define
\be  \om = e^{2\pi i/n}  \ee
and let $\ga$ be a primitive $n$th root of unity, i.e., a generator of the
group $\Ga_n$ (for example, $\ga = \om$). Let $\chi : \Ga_n \to \Co_{\ast}$
be the character defined by
\be  \la\chi,\ga^r\ra = \om^r\;\;\mbox{for all}\;\;r \in Z_n\,.  \ee
Then there exists a unique isomorphism of $\Ga_n$ onto $\Ga_{n}^{\ast}$
which maps $\ga$ onto $\chi$\,, and we could use this isomorphism to
identify $\Ga_{n}^{\ast}$ with $\Ga_n$\,. The canonical pairing is given by
\be  \la\chi^r,\ga^s\ra = \om^{rs}\;\;\mbox{for all}\;\;r,s \in Z_n\,, \ee
and an arbitrary bicharacter on $\Ga_n$ has the form
\be  \si_k(\ga^r,\ga^s) = \om^{krs}\;\;\mbox{for all}\;\;r,s \in Z_n\,, \ee
with some element $k \in Z_n$\,. According to Eq. \reff{ans}, the universal
$R$--matrix $R_k \in \Co\Ga_{n}^{\ast} \ot \Co\Ga_{n}^{\ast}$ corresponding
to $\si_k$ is given by
\be  R_k = \frac{1}{n^2}\sum_{p,q,r,s \in Z_n}\om^{kpq}\om^{-pr}\om^{-qs}
                                                  \chi^r \ot \chi^s\,.  \ee
Obviously, part of this sum can be carried out. Consider, for a moment, $k$
as an element of $\{1,2,\ldots,n\}$ and let $d$ be the greatest common
divisor of $k$ and $n$\,. Let $\ell$ be the inverse of $k/d$ in the ring
$Z_{n/d}$\,. Then
\be  R_k = \frac{d}{n}\sum_{x,\,y \in Z_{n/d}}(\om^d)^{-\ell xy}
                                (\chi^d)^x \ot (\chi^d)^y\,. \label{Rex} \ee
Note that in this case the subgroups $\De_1$ and $\De_2$ of $\Ga_{n}^{\ast}$
introduced in Proposition 2.2 are equal to
\be  \De_1 = \De_2 = \{(\chi^d)^x\,|\, x \in Z_{n/d}\}  \ee
and are isomorphic to $\Ga_{n/d}$\,. The universal $R$--matrix \reff{Rex}
with $k=1$ (and hence $d=1$, $\ell=1$) has already been used in Ref. [10],
the special case $n=2$ (related to supersymmetry) is contained in
Ref. [11]. \\

We would now like to interpret our results from a dual point of view.
However, this necessitates some preparation. Consequently, we shall first
comment on the duality of Hopf algebras and return to our main topic not
until Section 4.

\section{Duality of Hopf algebras and
     \newline
         coquasitriangular structures}
\setcounter{equation}{0}

The results of the present section are well--known to experts in Hopf algebra
theory. They are included for the convenience of the reader and for later
reference.

Let $A$ and $H$ be two bialgebras [2] and let
\be  \vp : A \times H \lra \Co  \ee
be a bilinear form on $A \times H$. Then there exists a unique bilinear
form
\be  \vp^{\ot2}: (A \ot A) \times (H \ot H) \lra \Co  \ee
such that
\be  \vp^{\ot2}(a \ot a',h \ot h') = \vp(a,h)\,\vp(a',h') \label{quad} \ee
for all $a,a' \in A$ and $h,h' \in H$.

The form $\vp$ is called a Hopf (or bialgebra) pairing of $A$ and $H$ if the
following relations are satisfied:
\bea  \vp^{\ot2}(\De_{A}(a),h \ot h') \eqq \vp(a,hh') \\
      \vp^{\ot2}(a \ot a',\De_{H}(h)) \eqq \vp(aa',h) \eea
\vspace{-6ex}
\bea  \vp(a,1_{H}) \eqq \ve_{A}(a) \\
      \vp(1_{A},h) \eqq \ve_{H}(h)  \eea
for all $a,a' \in A$ and all $h,h' \in H$. Here, $\De_A$, $\ve_A$, and
$1_A$ denote the coproduct, the counit, and the unit element of $A$\,, and
$\De_H$, $\ve_H$, and $1_H$ have the analogous meaning.

{}From the point of view of dual Hopf algebras (see Ref. [2] and below) the
concept of a Hopf pairing is very natural and has been used by many authors.
The first reference I know of is Ref. [12].

The prototype of Hopf pairings can be constructed as follows. Let $H$ be a
bialgebra, let $H^{\ast}$ be the vector space dual to $H$, and let
$H^{\circ} \subset H^{\ast}$ be the set of all linear forms on $H$ which
vanish on an ideal of $H$ of finite codimension. It is easy to see that
$H^{\circ}$ is a subspace of $H^{\ast}$. Ref. [2] contains various
characterizations of the elements of $H^{\circ}$. In particular,
let $m:H \ot H \rar H$ be the product mapping of $H$ and let
$m^{t}: H^{\ast} \rar (H \ot H)^{\ast}$ be its transpose. Recall that
$H^{\ast} \ot H^{\ast}$ is canonically embedded in $(H \ot H)^{\ast}$. Then
an element $f \in H^{\ast}$ belongs to $H^{\circ}$ if and only if
$m^t(f) \in H^{\ast} \ot H^{\ast}$.

Now let $\vp$ be the restriction of the canonical pairing
\be  \la\,\;,\;\ra:H^{\ast} \times H \lra \Co  \ee
onto $H^{\circ} \times H$. Then there exists a unique bialgebra structure
on $H^{\circ}$ such that $\vp$ is a Hopf pairing of $H^{\circ}$ and $H$ (for
a proof, see Ref. [2]). We call the bialgebra $H^{\circ}$ the Hopf dual
of $H$.

Using this result, the Hopf pairings can be characterized as follows. Let us
use the notation introduced at the beginning of this section. It is
well--known that, for any bilinear form $\vp$ on $A \times H$, there exist
two linear mappings
\be  \vp_{\ell}:A \lra H^{\ast} \ee
\be     \vp_{r}:H \lra A^{\ast} \ee
such that
\be  \la\vp_{\ell}(a),h\ra = \la\vp_{r}(h),a\ra = \vp(a,h)  \ee
for all $a \in A$ and $h \in H$ (obviously, they are uniquely determined by
this equation). Then we can prove the following proposition. \\

\noindent
{\bf Proposition 3.1} \\
{\it We use the notation introduced above. Then the following statements are
equivalent: \\
1) The bilinear form $\vp$ is a Hopf pairing of $A$ and $H$. \\
2) We have $\vp_{\ell}(A) \subset H^{\circ}$, and the map of $A$ into
$H^{\circ}$ induced by $\vp_{\ell}$ is a bialgebra homomorphism. \\
3) We have $\vp_{r}(H) \subset A^{\circ}$, and the map of $H$ into
$A^{\circ}$ induced by $\vp_{r}$ is a bialgebra homomorphism.}

Proof \\
The only non--trivial part of this proposition is that, for any Hopf pairing,
we have
\be  \vp_{\ell}(A) \subset H^{\circ}\;\;\;
                        \mbox{and}\;\;\;\vp_{r}(H) \subset A^{\circ}\,, \ee
and this follows from the characterization of $H^{\circ}$ and $A^{\circ}$
recalled above. The rest can be checked by straightforward calculations. \\

Suppose now that $H$ is even a Hopf algebra, and let $S_H$ be the antipode
of $H$. Then it is easy to see that the transpose $S_{H}^{t}$ of $S_H$ maps
$H^{\circ}$ into itself and that the map of $H^{\circ}$ into itself induced
by $S_{H}^{t}$ is an antipode for $H^{\circ}$. Thus if $H$ is a Hopf algebra,
then $H^{\circ}$ is likewise.

We now can prove the following corollary to Proposition 3.1. \\

\noindent
{\bf Corollary 3.1} \\
{\it We use the notation introduced at the beginning of this section, but
suppose in addition that $A$ and $H$ are Hopf algebras, with antipodes $S_A$
and $S_H$, respectively. If $\vp$ is a Hopf pairing of $A$ and $H$, then
\be  \vp(S_{A}(a),h) = \vp(a,S_{H}(h))  \ee
for all $a \in A$ and $h \in H$.}

Proof \\
It is well--known [2] that a bialgebra homomorphism of a Hopf algebra into
another one is, in fact, a Hopf algebra homomorphism, i.e., it intertwines
the antipodes. Thus the corollary follows from Proposition 3.1. \\

Our next task is to dualize the notion of a universal $R$--matrix [1]. Let
$H$ be a bialgebra. An element $R \in H \ot H$ is called a universal
$R$--matrix for $H$ if it has the following four properties:
\be  R\;\;\mbox{is invertible in}\;\;H \ot H  \label{URM1} \ee
\vspace{-5.5ex}
\be (T\De_{H}(h))R = R\,\De_{H}(h)\;\;\mbox{for all}\;\;h \in H \label{URM2}\ee
\vspace{-7ex}
\bea (\De_{H} \ot id_{H})(R) \eqq R_{13}R_{23}     \label{URM3} \\
     (id_{H} \ot \De_{H})(R) \eqq R_{13}R_{12}\,.  \label{URM4} \eea
Our notation is standard:
\be T:H \ot H \lra H \ot H  \ee
is the linear twist mapping given by
\be T(h \ot h') = h' \ot h\;\;\mbox{for all}\;\;h,h' \in H  \ee
and the elements
\be R_{12}\,,\:R_{23}\,,\:R_{13} \in H \ot H \ot H  \ee
are defined by
\be              R_{12} = R \ot 1_H  \ee
\be              R_{23} = 1_H \ot R  \ee
\be    R_{13} = (T \ot id_H)(R_{23}) = (id_H \ot T)(R_{12})\,,  \ee
where $1_H$ denotes the unit element of $H$.

In order to dualize these four properties, we assume that $H$ is the Hopf
dual of a bialgebra $A$\,,
\be    H = A^{\circ}\,.  \ee
Since $A$ and hence $A \ot A$ are coalgebras, their duals $A^{\ast}$ and
$(A \ot A)^{\ast}$ are algebras, and $A^{\circ}$ is a subalgebra of
$A^{\ast}$. Actually, we have the following chain of canonical injective
algebra homomorphisms:
\be A^{\circ} \ot A^{\circ} \lra A^{\ast} \ot A^{\ast}
                                   \lra (A \ot A)^{\ast}\,. \label{inj} \ee
Thus, since $R$ is an element of $A^{\circ} \ot A^{\circ}$, it can also be
considered as an element of $(A \ot A)^{\ast}$. On the other hand, due to the
familiar properties of tensor products of vector spaces, $(A \ot A)^{\ast}$
can be canonically identified with the space $B(A,A)$ of all bilinear forms
on $A \times A$\,. In particular, $B(A,A)$ inherits from $(A \ot A)^{\ast}$
the structure of an associative algebra with a unit element. The
multiplication is the so--called convolution. It is denoted by an asterisk
$\ast$ and is defined as follows. Let $\psi,\,\psi'$ be two elements of
$B(A,A)$ and let $a,b \in A$\,. If $\De_A$ is the coproduct of $A$ and if
\be  \De_{A}(a) = \sum_{i}a^1_i \ot a^2_i  \label{dela} \ee
\be  \De_{A}(b) = \sum_{j}b^1_j \ot b^2_j  \label{delb} \ee
with $a^1_i,a^2_i,b^1_j,b^2_j \in A$\,, we have
\be(\psi \ast \psi')(a,b) = \sum_{i,j}\psi(a^1_i,b^1_j)
                                   \,\psi'(a^2_i,b^2_j) \,. \label{conv} \ee
The unit element is the bilinear form $\ve_B$ on $A \times A$ defined by
\be \ve_B(a,b) = \ve_A(a)\ve_A(b) = \ve_A(ab)  \ee
for all $a,b \in A$\,, where $\ve_A$ denotes the counit of $A$\,.

Now if
\be R = \sum_{i}R^{1}_{i} \ot R^{2}_{i}  \ee
is an arbitrary element of $A^{\circ} \ot A^{\circ}$, the corresponding
bilinear form $\vr \in B(A,A)$ is given by
\be  \vr(a,b) = \sum_{i}\la R^{1}_{i},a\ra\la R^{2}_{i},b\ra  \ee
for all $a,b \in A$\,, where $\la\,\;,\;\ra$ denotes the canonical pairing of
$A^{\ast}$ and $A$\,.

If $R$ is invertible in $A^{\circ} \ot A^{\circ}$, then all the more in
$B(A,A)$. This means that there exists a bilinear form $\vr' \in B(A,A)$
such that
\be  \vr \ast \vr' = \vr' \ast \vr = \ve_B \,. \label{inv}  \ee
Note that if $A$ is finite--dimensional, then the homomorphisms in Eq.
\reff{inj} are bijective, and hence $R$ is invertible in
$A^{\circ} \ot A^{\circ}$ if and only if a bilinear form $\vr'$ of the type
described above exists. In general, however, we cannot expect that the
invertibility of $\vr$ in $B(A,A)$ implies the invertibility of $R$ in
$A^{\circ} \ot A^{\circ}$.

The remaining conditions \reff{URM2} -- \reff{URM4} can also be reformulated
in terms of certain convolution products involving $\vr$\,. Because of lack
of space, we prefer to present the results more explicitly, as follows.

Condition \reff{URM2} is satisfied if (and, if $A^{\circ}$ separates the
elements of $A$\,, only if) the following equation holds for all
$a,b \in A$\,:
\be  \sum_{i,j}b^{1}_{j} a^{1}_{i} \vr(a^{2}_{i},b^{2}_{j}) =
     \sum_{i,j}\vr(a^{1}_{i},b^{1}_{j}) a^{2}_{i} b^{2}_{j} \label{comm} \ee
(here and in the subsequent equations we are using the notation introduced
in Eqs. \reff{dela}, \reff{delb}). Moreover, Eq. \reff{URM3} is satisfied
if and only if
\be  \vr(bc,a) = \sum_{i}\vr(b,a^1_i)\vr(c,a^2_i) \label{rho3} \ee
for all $a,b,c \in A$\,, and Eq. \reff{URM4} holds if and only if
\be  \vr(a,bc) = \sum_{i}\vr(a^1_i,c)\vr(a^2_i,b) \label{rho4} \ee
for all $a,b,c \in A$\,.

The foregoing results motivate the following definition. \\

\noindent
{\bf Definition 3.1} \\
{\it Let $A$ be a bialgebra. A bilinear form $\vr$ on $A \times A$ is said
to define a coquasitriangular (or dual quasitriangular, or braiding)
structure on $A$ if it is convolution invertible (i.e., if it has an inverse
in the sense of Eq. \reff{inv}) and if the relations \reff{comm}, \reff{rho3},
and \reff{rho4} are satisfied. A bialgebra endowed with a coquasitriangular
structure will be called coquasitriangular.} \\

Apparently, dual quasitriangular structures have first been introduced in
Ref. [13] and shortly afterwards by various authors. In particular, I refer
the reader to Ref. [14] (note that the coquasitriangular structures discussed
here are the right braiding structures defined there).

Visibly, Eq. \reff{comm} is some sort of generalized commutation relation.
In fact, the famous RTT relations of Ref. [15] are of this type, and the
corresponding bi\-algebras $A(R)$ (with $R$ an invertible solution of the
Yang--Baxter equation) are coquasitriangular.

The Eqs. \reff{rho3} and \reff{rho4} can be rewritten as follows. Let
$\vr^{\ot 2}$ be the bilinear form on $A \ot A$ defined by
\be  \vr^{\ot 2}(a \ot a',b \ot b') = \vr(a,b)\vr(a',b')  \ee
for all $a,a',b,b' \in A$ (see Eq. \reff{quad}). Then Eq. \reff{rho3} is
equivalent to
\be  \vr^{\ot 2}(b \ot c, \De_{A}(a)) = \vr(bc,a)\,, \label{rho3mod} \ee
and Eq. \reff{rho4} is equivalent to
\be  \vr^{\ot 2}(T\De_{A}(a), b \ot c) = \vr(a, bc)\,, \label {rho4mod} \ee
(where, once again, $T$ denotes the twist mapping) or, what amounts to the
same, to
\be  \vr^{\ot 2}(\De_{A}(a), c \ot b) = \vr(a,bc)\,. \label{rho4modd} \ee

As is well--known, there is a real host of properties and relations satisfied
by a universal $R$--matrix. Thus one should expect the same for
coquasitriangular structures. Of course, due to the well--known problems
with duality for infinite--dimensional vector spaces, we cannot simply prove
these latter results by duality. Rather, we use duality to find out which
properties should be true, and then try to prove these directly.

Let us give some examples. In the following, $R$ denotes a universal
$R$--matrix for a bialgebra $H$ and $\vr$ denotes a coquasitriangular
structure for a bialgebra $A$\,. It is obvious that with $R$ also
$T(R^{-1}) = T(R)^{-1}$ is a universal $R$--matrix for $H$ (as before, $T$
denotes the appropriate twist mapping, here the one of $H \ot H$ onto
itself). Dually, with $\vr$ also $\vr' \comp T = (\vr \comp T)'$ is a
coquasitriangular structure for $A$\,, where the prime denotes the convolution
inverse. A coquasitriangular structure $\vr$ is said to be cotriangular if
$\vr' \comp T = \vr$\,.

It is known that
\be   (\ve_H \ot id_H)(R) = 1_H  \ee
\be   (id_H \ot \ve_H)(R) = 1_H  \ee
(in these equations, $\Co \ot H$ and $H \ot \Co$ have been canonically
identified with $H$). Dually, it can be shown that
\be   \vr(1_A,a) = \vr(a,1_A) = \ve_A(a) \label{rhounit} \ee
for all $a \in A$ (for example, apply Eq. \reff{rho4} with $c = 1_A$ and use
that $\vr$ has a convolution inverse).

Finally, let us assume that $H$ and $A$ are Hopf algebras. Then it is known
that
\bea  (S_H \ot id_H)(R) \eqq R^{-1}  \\
      (id_H \ot S_H)(R^{-1}) \eqq R\,.  \eea
Dually, it can be shown that
\bea  \vr(S_A(a),b) \eqq \vr'(a,b)  \label{rhoS1} \\
      \vr'(a,S_A(b)) \eqq \vr(a,b)  \eea
for all $a,b \in A$\,, where, as before, $\vr'$ is the convolution inverse
of $\vr$\,. In particular, we have

\be   \vr(S_A(a),S_A(b)) = \vr(a,b)  \ee
for all $a,b \in A$\,.

Conversely, if $\vr$ is a bilinear form on a Hopf algebra $A$ satisfying
the Eqs. \reff{comm} -- \reff{rho4} and \reff{rhounit}, then the bilinear
form $\vr'$ defined by Eq. \reff{rhoS1} is a convolution inverse of $\vr$
and hence $\vr$ defines a coquasitriangular structure on $A$\,.

Let us next make contact with our discussion of Hopf pairings. Let $A$ be a
bi\-algebra and let $A^{cop}$ (resp. $A^{aop}$) be the bialgebra which,
considered as an algebra (resp. a coalgebra) coincides with $A$\,, but which,
considered as a coalgebra (resp. an algebra) has the structure opposite to
that of $A$\,. Then the Eqs. \reff{rho3mod} -- \reff{rho4modd} and
\reff{rhounit} show that a coquasitriangular structure $\vr$ of $A$ can be
defined to be a Hopf pairing $A^{cop} \times A \rar \Co$ (or, equivalently, a
Hopf pairing $A \times A^{aop} \rar \Co$) which has a convolution inverse in
the sense of Eqs. \reff{conv}, \reff{inv} and satisfies the generalized
commutation relations \reff{comm}. In particular, it follows from
Proposition 3.1 that, for any $b \in A$\,, the linear forms $a \rar \vr(a,b)$
and $a \rar \vr(b,a)$ on $A$ belong to $A^{\circ}$.

This result is highly welcome. In our dualization process, we started from
a universal $R$--matrix $R \in A^{\circ} \ot A^{\circ}$ and embedded
$A^{\circ} \ot A^{\circ}$ into $(A \ot A)^{\ast} \simeq B(A,A)$. For
infinite--dimensional bialgebras, $(A \ot A)^{\ast}$ is terribly much larger
than $A^{\circ} \ot A^{\circ}$. On the other hand, a certain enlargement of
$A^{\circ} \ot A^{\circ}$ is reasonable, since for the most interesting
examples of infinite--dimensional quasitriangular Hopf algebras $H$, the
universal $R$--matrix does not belong to $H \ot H$ but to a certain
completion of it. Our last result shows that $\vr$ seems to be not as far
from $A^{\circ} \ot A^{\circ}$ as one might have feared. Thus one could hope
to find a subbialgebra $U$ of $A^{\circ}$ associated to $\vr$ and a suitable
completion of $U \ot U$ in $B(A,A)$ containing $\vr$\,. Such a programme has
been carried out in Ref. [14].

Finally, let us recall that, for any coalgebra $C$\,, the right $C$--comodules
may also be regarded as left $C^{\ast}$--modules (where $C^{\ast}$ denotes
the associative algebra dual to $C$)[2]. In view of the discussion of the
present section we expect that a coquasitriangular structure on a bialgebra
$A$ can be used to convert the class of right $A$--comodules into a braided
monoidal category. This can in fact be done (actually, this is the starting
point of Ref. [14]). We do not want to go into detail here but only show how
the corresponding braiding isomorphisms are defined.

Thus let $V$ and $W$ be two right $A$--comodules, with structure maps
\bea  \de_{V} \!\!\!& : &\!\!\! V \lra V \ot A \\
      \de_{W} \!\!\!& : &\!\!\! W \lra W \ot A \,. \eea
To define the braiding isomorphism
\be \psi_{V,W} : V \ot W \lra W \ot V \:, \ee
choose $x \in V$, $y \in W$, and write
\bea  \de_{V}(x) \eqq \sum_{i}x_i \ot a_i  \\
      \de_{W}(y) \eqq \sum_{j}y_j \ot b_j  \eea
with $x_i \in V$, $y_j \in W$\,; $a_{i},b_{j} \in A$\,. Then we set
\be  \psi_{V,W}(x \ot y) = \sum_{i,j}\vr(a_{i},b_{j})\, y_j \ot x_i \;.
                                                         \label{cobraid} \ee

We close this section by two remarks. For reasons of uniformity, we have
formulated the results of this section over the base field $\Co$\,. Of
course, everything holds for arbitrary base fields of any characteristic.
Moreover, by systematically using the techniques supplied by Ref. [3], these
results can easily be generalized to the graded case, with arbitrary Abelian
groups of degrees and arbitrary commutation factors.

\section{Coquasitriangular structures for the group
     \newline
         Hopf algebra of an Abelian group}
\setcounter{equation}{0}

We now are ready to rederive the results of Section 2 from the dual point of
view. Moreover, we are going to see that the dual picture is applicable also
to infinite Abelian groups.

To begin with we still assume that the group $\Ga$ is finite. Then the group
Hopf algebras $\CG$ and $\CGS$ are (Hopf) dual to each other. In fact, since
$\Ga$ is a basis of $\CG$ and $\Gas$ a basis of $\CGS$, there is a unique
bilinear form $\CGS \times \CG \rar \Co$ whose restriction to
$\Gas \times \Ga$ is equal to the group pairing $\Gas \times \Ga \rar \Co$\,,
and it is easy to check that this bilinear form is a non--degenerate Hopf
pairing. Thus there is no problem in denoting these two pairings by the same
symbol $\la\,\;,\;\ra$.

The result above is well--known, it agrees with the fact that the dual of
the group Hopf algebra $\Co G$ of a finite group $G$ is canonically isomorphic
to the Hopf algebra ${\cal F}(G)$ of complex--valued functions on $G$\,. In
fact, the characters of $\Ga$ form a basis of ${\cal F}(\Ga)$. Of course, all
this is closely related to Fourier transformation on $\Ga$ and $\Gas$.

According to Section 3 the foregoing result implies that the universal
$R$--matrices for $\CGS$ are in bijective correspondence with the
coquasitriangular structures on $\CG$. More precisely, consider the chain
\reff{inj} of injective algebra homomorphisms, with $A = \CG$. In the present
case these are, in fact, algebra isomorphisms, and according to the foregoing
result we can (and will) identify $(\CG)^{\ast}$ with $\CGS$. On the other
hand $(\CG \ot \CG)^{\ast}$ is canonically isomorphic to $B(\CG,\CG)$, the
algebra of bilinear forms on $\CG \times \CG$. Since $\Ga$ is a basis of
$\CG$, the bilinear forms on $\CG \times \CG$ are uniquely determined by
their values on $\Ga \times \Ga$, which shows that the vector space
$B(\CG,\CG)$ can be identified with ${\cal F}(\Ga \times \Ga)$, the vector
space of complex--valued functions on $\Ga \times \Ga$. Summarizing, we have
the following chain of vector space isomorphisms:
\be \CGS \ot \CGS \lra (\CG \ot \CG)^{\ast} \lra B(\CG,\CG) \lra
                                             {\cal F}(\Ga \times \Ga)\,. \ee
We know that the first of these maps is an algebra homomorphism and (by
definition of the convolution) so is the second. A look at Eq. \reff{conv}
shows that the third map is an algebra homomorphism as well. Furthermore, it
is easy to check that the inverse mapping
\be  {\cal F}(\Ga \times \Ga) \lra \CGS \ot \CGS  \ee
is just the mapping $\si \rar R_{\si}$ mentioned in part 1) of
Proposition 2.1.

Thus all we have to do is to determine those bilinear forms $\vr$ on
$\CG \times \CG$ which define a coquasitriangular structure on $\CG$. Since
$\Ga$ is Abelian, the condition \reff{comm} is trivially satisfied.
Obviously, the remaining conditions (convolution invertibility and Eqs.
\reff{rho3}, \reff{rho4}) are fulfilled if and only if the restriction of
$\vr$ to $\Ga \times \Ga$ (i.e., the element of ${\cal F}(\Ga \times \Ga)$
corresponding to $\vr$) is a bicharacter on $\Ga$. Thus we have rederived
part 1) of Proposition 2.1 and the non--trivial parts of Corollary 2.1.

{\it Let us now drop the assumption that $\Ga$ is finite.} Then the
coquasitriangular structures on $\CG$ can be characterized as before: A
bilinear form $\vr$ on $\CG \times \CG$ defines a coquasitriangular structure
on $\CG$ if and only if its restriction $\si$ to $\Ga \times \Ga$ is a
bicharacter of $\Ga$, and $\vr$ is cotriangular if and only if $\si$ is a
commutation factor on $\Ga$.

Of course, we should now answer the question of how graded multilinear
algebra (which has been the point of departure for our round trip) can be
formulated in the present setting. The answer is suggested by the remarks
at the end of Section 3: The $\CG$--comodules should take the role of the
$\CGS$--modules.

In fact, it is well--known that $\Ga$--graded vector spaces and
$\CG$--comodules are essentially the same thing. Indeed, let $V$ be a right
$\CG$--comodule and let
\be   \de : V \lra V \ot \CG  \ee
be its structure mapping. If we define, for all $\ga \in \Ga$,
\be   V_{\ga} = \{x \in V \,|\, \de(x) = x \ot \ga \}\,,  \ee
then, obviously, $V_{\ga}$ is a subspace of $V$, and it is easy to see that
these subspaces form a $\Ga$--gradation of $V$\,:
\be   V = \bigoplus_{\ga \in \Ga} V_{\ga}\,.  \ee
Conversely, if $(V_{\ga})_{\ga \in \Ga}$ is a $\Ga$--gradation of a vector
space $V$, there exists a unique linear mapping
\be   \de : V \lra V \ot \CG  \ee
such that
\be   \de(x) = x \ot \xi  \ee
for all $x \in V_{\xi}$\,, $\xi \in \Ga$, and this mapping defines on $V$
the structure of a right $\CG$--comodule. Obviously, these transitions from
right $\CG$--comodules to $\Ga$--graded vector spaces and vice versa are
inverse to each other. Moreover, the homomorphisms of $\CG$--comodules are
exactly the homomorphisms of $\Ga$--graded vector spaces.

Now let $\vr$ define a coquasitriangular structure on $\CG$ and let $\si$
be the restriction of $\vr$ to $\Ga \times \Ga$. Consider two $\Ga$--graded
vector spaces, i.e., right $\CG$--comodules, denoted by $V$ and $W$, say.
Then the braiding isomorphism
\be   \psi_{V,W} : V \ot W \lra W \ot V  \ee
corresponding to $\vr$ (see Eq. \reff{cobraid}) is given by
\be   \psi_{V,W}(x \ot y) = \si(\xi,\eta)\,y \ot x  \ee
for all $x \in V_{\xi}$\,, $y \in W_{\eta}$\,; $\xi,\eta \in \Ga$. This is
exactly the formula familiar from graded multilinear algebra [3].

Needless to say, there is a lot of other details which should be checked
before we can be sure that the graded multilinear algebra as developed in
Ref. [3] and the comodule picture described here (with $\vr$ cotriangular)
are really equivalent. However, I do not want to embark on this boring
exercise. Instead, I would like to close this section by the remark that we
have come across one more example of a well--known fact: In the
infinite--dimensional case, coalgebras and comodules are easier to handle
than algebras and modules.

\section{Discussion}
\setcounter{equation}{0}
In the present work we have described how the graded multilinear algebra as
developed in Ref. [3] can be understood in the language of triangular or
cotriangular Hopf algebras. We may now proceed to apply the theory of Hopf
algebras to the study of algebras living in the graded multilinear category.
In particular, we could apply Majid's theory of bosonization [16] and try
to reduce (or, at least, to relate) the theory of Lie colour algebras to the
theory of usual Lie algebras. (Roughly speaking, bosonization is a procedure
which converts a Hopf algebra living in a category with non--trivial braiding
into a usual Hopf algebra.) The idea to try this is not new. It is used in
supersymmetry by introducing Grassmann variables, and it is also applied in
Ref. [5], where the present author has reduced the theory of Lie colour
algebras to the theory of graded Lie (super)algebras.

More recent and closer to the present work and to the bosonization procedure
are the Refs. [17,\,9]. In Ref. [17] the authors prove, in the context of
coquasitriangular Hopf algebras, a general bosonization theorem and use it
to derive Schur's double centralizer theorem for cotriangular Hopf algebras
(hence, in particular, for Lie colour algebras). In Ref. [9] the author
constructs a bosonized version of the enveloping algebra of a Lie colour
algebra and determines its universal $R$--matrix.

Of course, there is no need to restrict attention to the triangular case,
i.e., to commutation factors. In fact, allowing for arbitrary bicharacters,
the graded vector spaces are the objects of a really braided monoidal
category, i.e., one in which the braid groups really come into play. On the
other hand, these categories are still very close to the classical case,
hence they might be used to become acquainted with the braided situation.
For example, one could try to define generalized Lie algebras living in these
categories. This has in fact been done in a recent work by Pareigis [18].
Using the tensor operator techniques developed in Refs. [19,\,20] one might
then hope to construct interesting quantum spin chain Hamiltonians which
are invariant under these algebras.\\[2ex]

\noindent
{\bf Acknowledgements}\,\,
A major part of the present work has been carried out during two visits of
the author to the Erwin Schr\"{o}dinger Institut in Vienna. The kind
invitations by Harald Grosse and the hospitality extended to the author in
the ESI are gratefully acknowledged. \\

\newpage

\noindent
{\large\bf References}
\begin{enumerate}
\item V. Drinfeld, J. Sov. Math. {\bf 41}, 898 (1988) (expanded version of a
      report to the International Congress of Mathematicians, Berkeley 1986).
\item M.E. Sweedler, {\it Hopf algebras}, W.A. Benjamin, New York (1969).
\item M. Scheunert, J. Math. Phys. {\bf 24}, 2658 (1983).
\item V. Rittenberg and D. Wyler, Nucl. Phys. B {\bf 139}, 189 (1978).
\item M. Scheunert, J. Math. Phys. {\bf 20}, 712 (1979).
\item S. Mac Lane, {\it Categories for the working mathematician}, Springer,
      New York, Heidelberg, Berlin (1971).
\item A. Joyal and R. Street, Braided monoidal categories, Mathematics
      Reports 860081, Macquarie University (1986), and Adv. in Math.
      {\bf 102}, 20 (1993).
\item S. Majid, Int. J. Mod. Phys. A {\bf 5}, 1 (1990).
\item D.S. McAnally, Lett. Math. Phys. {\bf 33}, 249 (1995).
\item S. Majid, in {\it Spinors, twistors, Clifford algebras, and quantum
      deformations,} Proceedings of the 2\,nd Max Born Symposium, Wroclaw
      1992, edited by Z. Oziewicz et al., Kluwer, pp. 327--336.
\item D. Radford, J. Algebra {\bf 141}, 354 (1991).
\item M. Takeuchi, {\it The \#--product of group sheaf extensions applied
      to Long's theory of dimodule algebras,} Algebra--Berichte,
      Mathematisches Institut der Universit\"{a}t M\"{u}nchen, Nr. 34,
      Verlag Uni--Druck, M\"{u}nchen (1977).
\item S. Majid, in {\it Quantum probability and related topics VI,} Trento
      1989, proceedings, edited by L. Accardi et al., World Scientific,
      Singapore (1991), pp. 333--358.
\item R.G. Larson and J. Towber, Comm. Algebra {\bf 19}, 3295 (1991).
\item N.Yu. Reshetikhin, L.A. Takhtadzhyan, and L.D. Faddeev, Leningrad
      Math. J. {\bf 1}, 193 (1990).
\item S. Majid, J. Algebra, {\bf 163}, 165 (1994).
\item D. Fischman and S. Montgomery, J. Algebra {\bf 168}, 594 (1994).
\item B. Pareigis, On Lie algebras in braided categories, Mathematisches
      Institut der Universit\"{a}t M\"{u}nchen, preprint, April 1995.
\item M. Scheunert, in {\it Generalized symmetries in physics,} Clausthal
      1993, proceedings, edited by H.--D. Doebner, V.K. Dobrev, and A.G.
      Ushveridze, World Scientific (1994), pp. 77--89.
\item M. Scheunert, Int. J. Mod. Phys. B {\bf 8}, 3655 (1994).
\end{enumerate}

\end{document}